\documentclass{amsart}
\usepackage{url}
\usepackage{harvard}

\DeclareMathOperator{\hol}{hol}

\DeclareMathOperator{\ev}{ev}
\DeclareMathOperator{\spec}{spec}
\DeclareMathOperator{\Man}{\mathbf {Man}}
\theoremstyle{plain}
\newtheorem{theorem}{Theorem}[section]

\newtheorem{proposition}[theorem]{Proposition}

\theoremstyle{definition}
\newtheorem{definition}[theorem]{Definition}

\theoremstyle{remark}
\newtheorem{example}{Example}[section]

\numberwithin{equation}{section}
\numberwithin{figure}{section}

\newcommand{\cH}{{\mathcal H}}
\newcommand{\cG}{{\mathcal G}}

\newcommand{\cA}{{\mathcal A}}

\newcommand{\cU}{{\mathcal U}}

\newcommand{\cP}{{\mathcal P}}

\newcommand{\cV}{{\mathcal V}}

\newcommand{\CC}{{\mathbb C}}
\newcommand{\RR}{{\mathbb R}}
\newcommand{\ZZ}{{\mathbb Z}}

\renewcommand{\a}{\alpha}
\renewcommand{\b}{\beta}
\renewcommand{\c}{\gamma}

\citationstyle{dcu}

\begin{document}

\title[Introduction]{An introduction to bundle gerbes}
\ \author[M.K. Murray]{Michael K Murray}
  \address[Michael K. Murray]
  {School of Mathematical Sciences\\
  University of Adelaide\\
  Adelaide, SA 5005 \\
  Australia}
  \email{michael.murray@adelaide.edu.au}
\thanks{The author acknowledges the support of the Australian
Research Council and would like to thank the organisers of the conference for the
invitation to participate. The referee is also thanked for useful remarks on the manuscript.}


\begin{abstract}
An introduction to the theory of bundle gerbes and their relationship to Hitchin-Chatterjee gerbes
is presented.  Topics covered are connective structures, triviality and stable isomorphism
as well as examples and applications. 
\end{abstract}
\maketitle

\section{Dedication}
Over the years I have had many interesting mathematical conversations with Nigel and regularly
came away with a solution to a problem or a new idea. While preparing this article I was trying to recall
when he first told me about gerbes.  I thought for awhile that age was going to get the better of my memory 
as many conversations seemed to have blurred together. But then  I discovered that the annual departmental
 research reports really do have their uses. In July of 1992 I attended the 
`Symposium on gauge theories and     topology'  at Warwick and reported in the 1992 Departmental Research Report that:           
\begin{quote}
 I \dots had discussions
with Nigel Hitchin about `gerbes'. These are a generalisation
of line bundles \dots
\end{quote}
Further searching of my electronic files revealed an order for Brylinski's book {\sl Loop spaces, characteristic classes 
and geometric quantization} on the 29th April 1993. I recall that the book took some months to makes its away across the
 sea to Australia during which time I pondered the advertising material I had which said that gerbes were fibrations of
  groupoids. Trying to interpret this  lead to a paper on  bundle gerbes which I submitted to Nigel in his role as a London 
  Mathematical Society Editor on the  25th July 1994. The Departmental Research Report of the same year reports that:
\begin{quote}  This year
I began some work on a geometric construction 
called a bundle gerbe. These provide a 
geometric realisation of the three dimensional 
cohomology of a manifold. 
\end{quote}
My sincere thanks to Nigel for introducing me to gerbes and for the many other fascinating insights into mathematics 
that he has given me over the  years.
         
\section{Introduction}
\label{sec:intro}

The theory of gerbes began with Giraud \citeyear{Gir} and was popularised in the book  by Brylinski \citeyear{Bry}. A short
 introduction by Nigel Hitchin \citeyear{Hit2003} in the `What is a \dots ?' series can be found in  the Notices of the AMS. 
 Gerbes provide a geometric realisation of the three dimensional cohomology of a manifold in a manner analogous to the way a 
 line bundle is a geometric realisation of two dimensional cohomology. Part of the reason for their recent popularity is applications 
 to string theory in particular the notion of the $B$-field. Strings on a manifold are elements in the loop space of the manifold and 
 we would expect their quantization to involve a hermitian line bundle on the loop space arising from a two class on the loop space. 
 That two class can arise as the transgression of some three class on the underlying manifold. Gerbes provide a geometrisation of this 
 process.  String theory however is not the only application of gerbes and we refer the interested reader to the related work of Hitchin \citeyear{MR1876068,MR2217300} which applies  gerbes to generalised geometry   and to reviews such as \cite{CarMicMur00Bundle-gerbes-applied} 
 and \cite{Mic2} which give applications of gerbes to other
problems in quantum field theory. 

As with everything else in the theory of gerbes, the  relationship of bundle gerbes to gerbes is best understood by comparison 
with the case of hermitian line bundles or equivalently $U(1)$ (principal) bundles. 
There are basically three ways of thinking about $U(1)$  bundles over a manifold $M$:
\begin{enumerate}
\item  A certain kind of locally free sheaf on $M$.
\item A co-cycle $g_{\a\b} \colon U_\a \cap U_\b \to U(1)$ for 
some open cover $\cU = \{ U_\a \mid \a \in I \}$ of $M$.
\item A principal $U(1)$ bundle $P \to M$.
\end{enumerate}
In the case of gerbes over $M$ we can think of these as:
\begin{enumerate}
\item  A certain kind of sheaf of groupoids on $M$ (Giraud, Brylinski). 
\item A co-cycle $g_{\a\b\c} \colon U_\a \cap U_\b \cap U_\c \to U(1)$ for 
some open cover $\cU = \{ U_\a \mid \a \in I \}$ of $M$ or alternatively a choice
of $U(1)$ bundle $P_{\a\b} \to U_\a \cap U_\b$ for each double overlap  (Hitchin, Chatterjee).
\item A bundle gerbe (Murray).
\end{enumerate}

Note that we are slightly abusing the definition of gerbe here as what we are considering are gerbes with 
band the sheaf of smooth functions from $M$ into $U(1)$. There are more general kinds of gerbes on $M$ just as there are more general
 kinds of sheaves on $M$ beyond those arising as the sheave of sections of a hermitian line bundle.  

Recall some of the basic facts about $U(1)$ bundles on a manifold $M$.  
\begin{enumerate}
\item If $P \to M$ is a $U(1)$ bundle there is a dual  bundle $P^* \to M$ and if $Q \to M$ is another
$U(1)$ bundle there is a {\em product} $P \otimes Q \to M$.
\item If $f \colon N \to M$ is a smooth map there is a {\em pullback} bundle $f^*(P) \to N$ and this behaves
well with respect to dual and product. That is $f^*(P^*) $ and $(f^*(P))^*$ are isomorphic as also are 
$f^*(P\otimes Q)$ and $f^*(P) \otimes f^*(Q)$. 
\item Associated to a $U(1)$ bundle $P \to M$ is a characteristic class, the chern class, $ c(P) \in H^2(M, \ZZ)$, which
is natural with respect to pullback, that is $f^*(c(P)) = c(f^*(P))$ and additive with respect to product and 
dual, that is $c(P \otimes Q) = c(P) + c(Q) $ and $c(P^*) = -c(P)$. 
\item $P \to M$ is called {\em trivial} if it is isomorphic to $M \times U(1)$ or equivalently admits a global section. 
$P$ is trivial if and only if $c(P) = 0$. 
\item There is a notion of a {\em connection}  on $P \to M$. Associated to a connection $A$ on $P$ is a closed two-form 
$F_A$ called the {\em curvature} of $A$ with the property that $F_A/2\pi i$ is a de Rham representative for the image
of $c(P)$ in real cohomology. 
\item If $\gamma \colon S^1 \to M$ is a loop in $M$ and $P \to M$ a $U(1)$  bundle with connection $A$ then parallel transport
 around $\gamma$ defines the 
{\em holonomy}, $\hol(A, \gamma)$ of $A$ around $\gamma$ which is an element of $U(1)$. If $\gamma$ is the boundary of 
a disk $D\subset M$ then we have
$$
\hol(A, \partial D) = \exp\left(\int_D F_A\right).
$$
\end{enumerate}

A gerbe is an attempt to generalise all the above facts about $U(1)$ bundles to some new kind of mathematical 
object in such a way that the characteristic class is in three dimension  cohomology.  Obviously for consistency other
 dimensions then 
have to  change. In particular the curvature should be a three-form and holonomy should be over two dimensional submanifolds. 
It turns out to be useful to consider the 
general case of any dimension of cohomology which we call a $p$-gerbe.  For historical reasons a $p$-gerbe has a 
characteristic class in $H^{p+2}(M, \ZZ)$ so the interesting  values of $p$ are $-2, -1, 0, 1, \dots$ with $U(1)$
bundles corresponding to $p=0$. 

A $p$-gerbe
then is some mathematical object which represents $p+2$ dimensional cohomology.  To make  completely precise
 what representing $p+2$ dimensional cohomology means  would take us to far afield from the present topic but we give a sketch 
here to motivate the behaviour we are looking for in $p$-gerbes.    To this end we will assume our $p$-gerbes $P$ live in some category
$\cG$ and  there is a (forgetful) functor $\Pi \colon \cG \to \Man$ the category of manifolds.   The functor
$\Pi$ and the category $\cG$ have to satisfy:

\begin{enumerate}
\item If $P$ is a $p$-gerbe there is a dual  $p$-gerbe $P^*$ and if $Q $ is another
$p$-gerbe there is a {\em product} $p$-gerbe $P \otimes Q $. In other words $\cG$ is monoidal and 
has a dual operation.
\item If $f \colon N \to M$ is a smooth map and $\Pi(P) = M$ there is a {\em pullback} $p$-gerbe $f^*(P) $
and a morphism $\hat f \colon f^*(P) \to P$ such that  $\Pi(f^*(P)) = N$ and $\pi(\hat f) = f$.  Pullback should 
behave well  with respect to dual and product. That is $f^*(P^*) $ and $(f^*(P)^*$ should be isomorphic as also should be
$f^*(P\otimes Q)$ and $f^*(P) \otimes f^*(Q)$. 
\item Associated to a $p$-gerbe $P$ is a characteristic class $ c(P) \in H^{p+2}(\Pi(P), \ZZ)$, which
is natural with respect to pullback, that is $f^*(c(P)) = c(f^*(P))$ and additive with respect to product and 
dual, that is $c(P \otimes Q) = c(P) + c(Q) $ and $c(P^*) = -c(P)$. 
\item As well as the notion of $P$ and $Q$ being isomorphic there is a possibly weaker notion of {\em equivalence}
where $P$ and $Q$ are equivalent if and only if $c(P) = c(Q)$.  We say $P$ is {\em trivial} if $c(P) = 0$. 
\item There is a notion of a {\em connective structure}  $A$ on $P$. Associated to a connective structure $A$ on $P$ is a closed $(p+2)$-form 
$\omega$ on $\Pi(P)$ called the {\em ($p+2$)-curvature} of $A$ with the property that $\omega/2\pi i$ is a de Rham representative
 for the image of $c(P)$ in real cohomology. 
\item If $X \subset  \Pi(P)$ is an oriented $p+1$ dimensional submanifold of  $\Pi(P)$ we should be able to  define the
holonomy of the connective structure  $\hol(A, X) \in U(1)$ over $X$.  Moreover if $  Y \subset  \Pi(M)$ is  an oriented 
$p+2$-dimensional submanifold with boundary then we want to have that
$$
\hol(\cA, \partial Y) = \exp\left(\int_Y F_\omega\right).
$$
\end{enumerate}

Clearly by construction the category of $U(1)$ bundles, with the forgetful functor which  assigns to a $U(1)$ bundle
its base manifold, is an example of a $0$-gerbe. 
Before we consider other examples we need some facts about bundles with structure group an abelian Lie group $H$. 
If $P \to M$ is an $H$ bundle on $M$ then by choosing local sections of $P$ for an open cover $\cU = \{ U_\a \mid \a \in I \}$ we
 can construct transition  functions $h_{\a\b} \colon U_\a \cap U_\b  \to H$ and in the usual way this defines a 
class  $c(P) \in H^1(M, H)$ where here we  abuse notation and write $H$ for what is really the sheaf of smooth functions with
 values in $H$. It is a standard fact that 
isomorphism classes of $H$ bundles are in bijective correspondence with $H^1(M, H)$ in this manner.  If $P \to M$ is an $H$ 
bundle we can define its {\em dual} as follows. Let $P^*$ be isomorphic to $P$ as a manifold with projection to $M$ and 
for convenience let $p^* \in P^*$ denote $p \in P$ thought of as an element of $P^*$. Then define a new $H$ action 
on $P^*$ by $p^* h = (ph^{-1})^*$. It is ovious that if $h_{\a\b}$ are transition functions for $P$ then $h_{\a\b}^* = h_{\a\b}^{-1}$ 
are transition functions for $P^*$. In particular we have that $c(P^*) = -c(P)$ if we write the group structure on $H^1(M, H)$ additively. 
If $Q$ is another $H$ bundle we can form the fibre product $P \times_M Q$ and let $H$ act on it by 
$(p, q)h = (ph, qh^{-1})$. Denote the orbit of $(p, q)$ under this action by $[p, q]$ and define an $H$ action 
by $[p, q]h = [p, qh] = [ph, q]$. The resulting $H$ bundle is denoted by $P \otimes Q \to M$.  If $h_{\a\b}$ are transition 
functions for $P$ and $k_{\a\b}$ are transition functions for $Q$ then $h_{\a\b} k_{\a\b}$ are
 transition functions for $P \otimes Q$ and thus $c(P \otimes Q ) = c(P) + c(Q)$. 
Notice that these constructions will not generally work for non-abelian groups because in such a case the action of $H$ on $P^*$ is 
not a right action and the action on $P \otimes Q$ is not even well-defined. 

\begin{example} 
The simplest example is that of functions $f \colon M \to \ZZ$ where we define the  functor $\Pi$ by $\Pi(f) = M$. 
The degree of $f$ is the class induced in $H^0(M, \ZZ)$ so functions from $M$ to $\ZZ$ are $-2$ gerbes over $M$. 
 Product and dual are pointwise addition and negation. There is no sensible
notion of connective structure.
\end{example}

\begin{example} Consider next  principal $\ZZ$ bundles $P \to M$. Clearly we want the functor $\Pi$ to be $\Pi(P) = M$
 and pull-backs are well known to exist.  As $\ZZ$ is abelian the constructions above apply and there are duals and 
 products.  The isomorphism class of a bundle is determined by a  class in $H^1(M, \ZZ)$ so $\ZZ$ bundles are $p=-1$ 
 gerbes. A $\ZZ$ bundle is trivial as a $-1$ gerbe if and only if it is trivial as a $\ZZ$ bundle. 

It is not immediately obvious what a connective structure on a $\ZZ$ bundle is but it  turns out that the correct notion
 is that of a $\ZZ$ equivariant map $\hat\phi \colon P \to i\RR$ where the action of $n \in \ZZ$ on $i\RR$ is addition
  of $2 \pi i n$ so that $\hat\phi(pn) = \hat\phi(p) + 2\pi i n$.  The map $\hat \phi$ then descends to a map
   $\phi \colon M \to S^1$ and the class of the bundle is the degree of this map.
The pull-back of the standard one-form on $\RR$, that is $d\hat\phi$ 
is a one-form on $P$ which descends to a one-form $\phi^{-1}d\phi$ on $M$. The de Rham class $(\phi^{-1}d\phi)/2\pi i$ is the image of 
the class of the bundle in real cohomology.  

We expect holonomy to be over a $-1 + 1 = 0$ dimensional submanifold. If $m_1, \dots, m_r$
is a collection of points in $M$ with each $m_i$ oriented by some $\epsilon_i \in \{ \pm 1\}$ let us 
denote by $\sum \epsilon_i m_i$ their union as an oriented zero dimensional submanifold of $M$. Then we define
$$
\hol\left( \hat\phi, \sum \epsilon_i m_i\right) = \prod_{i=1}^r \phi(m_i)^{\epsilon_i}.
$$
In the case of an oriented one-dimensional submanifold $X \subset M$ with ends $-X_0$ and $+X_1$  the fundamental theorem of calculus tells us that 
$$
\hol\left(\hat\phi, X_1 - X_0 \right) = \exp\left( \int_X d\phi \right)
$$ 

Notice that if we express a $\ZZ$ bundle locally in terms of transitions functions these are maps of the form
$f_{\a\b} \colon U_\a \cap  U_\b \to \ZZ $. That is, over  each double overlap we have a $-2$ gerbe. 
\end{example}

\begin{example}
It is clear from the above example that maps $\phi\colon M \to U(1)$ are also $-1$ gerbes with a connective structure.  The dual and product are just 
pointwise inverse and pointwise product. The class is the degree and the  connective structure is included automatically as part of $\phi$. 

We can also forget that there is a natural connective structure and just regard maps $\phi \colon M \to U(1)$ as $-1$ gerbes. 
In that case  the natural notion of isomorphism between two maps $\phi, \chi \colon M \to U(1)$ would be equality.
 However two such maps have the same degree if and only if they are homotopic. So the notion of equivalence of
  maps $\phi \colon M \to U(1)$, thought of as  $-1$ gerbes (without connective structure) should be homotopy and
   is different to the notion of isomorphism. 
\end{example}

\begin{example}
As we have remarked $U(1)$ bundles are, of course, $p=0$ gerbes. Notice that locally a $U(1)$ bundle is given by 
transition functions $g_{\a\b} \colon U_\a \cap  U_\b \to U(1) $, that is on each double overlap we have a $-1$ gerbe (with connective structure).
\end{example} 

We will see below that this pattern of a $p$ gerbe being defined as a $p-1$ gerbe on double overlaps of some open 
cover is exploited by Hitchin and Chatterjee to give  a definition of a $1$ gerbe.  But first we need some additional background material. 

\section{Background}

We will be interested in surjective submersions $\pi \colon Y \to M$ which we regard as generalizations of open covers. In particular 
if $\cU = \{ U_\a \mid \a \in I \}$ is an open cover we have the disjoint union 
$$
Y_{\cU} = \{(x, \a) \mid x \in U_\a \} \subset M \times I 
$$
with projection map $\pi(x, \a) = x$. The surjective morphism $\pi\colon Y_\cU \to M$ is called the {\em nerve} of the open cover $\cU$. 

A morphism of surjective submersions $\pi \colon Y \to M$ and $p \colon X \to M$ is a map 
$\rho \colon Y \to X$ covering the identity, that is $p \circ \rho = \pi$.  Any surjective submersion
 $\pi \colon Y \to M$ admits local sections so there is an open cover $\cU$ of $M$ and local sections
  $s_\a \colon U_\a \to Y$ of $\pi$. These local sections define a morphism $s \colon Y_\cU \to Y$ by $s(x, \a) = s_\a(x)$. 
  Indeed any morphism $Y_\cU \to Y$ will be of this form.    
If  $\cV = \{ V_\a \mid \a \in J \}$ is a refinement of $\cU$, that is there is  a map $\rho
\colon J \to I$ such that for every $\a \in J$ we have $V_{\a} \subset U_\rho({\a})$, we have  morphism of surjective
submersions $Y_\cV \to Y_\cU$ defined by $(\a, x) \mapsto (\rho(\a), x)$. 

Given a surjective morphism $\pi \colon Y \to M$ we can form the $p$-fold fibre product
$$
Y^{[p]} = \{(y_1, \dots, y_p) \mid \pi(y_1) = \cdots = \pi(y_p) \} \subset Y^p.
$$
The submersion property of $\pi$ implies that $Y^{[p]}$ is a submanifold of $Y^p$. There are smooth maps
$\pi_i \colon Y^{[p]} \to Y^{[p-1]}$, for $i=1, \dots, p$,  defined by omitting the $i$th element. 
 We will be interested in two particular examples.
 
 \begin{example} 
 If $\cU$ is an open cover of $M$ then the $p$-th fibre product $Y_\cU^{[p]}$ is the disjoint union of all the {\em ordered} 
 $p$-fold intersections. For example if $\cU = \{U_1, U_2\}$ is an open cover of $M$ then $Y_{\cU}^{[2]}$ is the disjoint
 union of $U_1 \cap U_2$, $U_2 \cap U_1$, $U_1 \cap U_1$ and $U_2 \cap U_2$. 
 \end{example}
 
 \begin{example} 
 If $P \to M$ is a principal $G$ bundle then $P \to M$ is a surjective submersion. It is easy to show that $P^{[p]} = P \times G^{p-1}$.
In particular $P^{[2]} = P \times G$ and we shall need later the related fact that there is a map $g \colon P^{[2]} \to G $ defined
by $p_1 g(p_1, p_2) = p_2$. 
 \end{example}
 
 Let $\Omega^q(Y^{[p]})$ be the space of differential $q$-forms on $Y^{[p]}$.  Define 
$$
\delta \colon \Omega^q(Y^{[p-1]}) \to \Omega^q(Y^{[p]}) \quad\text{by}\quad \delta =  \sum_{i=1}^p (-1)^{p-1} \pi_i^* .
$$

These maps form the {\em fundamental  complex}
$$
0 \to \Omega^q(M) \stackrel{\pi^*}{\to} \Omega^q(Y) \stackrel{\delta}{\to}  \Omega^q(Y^{[2]}) \stackrel{\delta}{\to} 
\Omega^q(Y^{[3]})    \stackrel{\delta}{\to} \dots 
$$
and from \cite{Mur96Bundle-gerbes} we have:

\begin{proposition}
The fundamental complex is exact for all $q \geq 0$. 
\end{proposition}

Note that if $Y = Y_\cU$ then this Proposition is a well-known result about the \v Cech de Rham double complex. 
See, for example Bott and Tu's book \citeyear{BotTu}.

Finally we need some notation. Let $H$ be an abelian group. 
If $g \colon Y^{[p-1]} \to H$ we define   $\delta(g) \colon Y^{[p]} \to H$ by 
$$
\delta(g) = (g \circ \pi_1) -  (g \circ \pi_2) +  ( g \circ \pi_3 ) \cdots .
$$

If $P \to Y^{[p-1]}$ is an $H$ bundle we define an $H$ bundle
$\delta(P) \to Y^{[p]}$ by
$$
\delta(P) = \pi_1^*(P) \otimes (\pi_2^*(P))^* \otimes \pi_3^*(P) \otimes \cdots .
$$
It is easy to check that $\delta(\delta(g)) = 1$ and  that $\delta(\delta(P)) $ is canonically trivial.

\section{Bundle gerbes}

\begin{definition} A {\em bundle gerbe} \cite{Mur96Bundle-gerbes} over $M$ is a pair $(P, Y)$ where $Y \to M$ 
is a surjective submersion and $P \to Y^{[2]}$ is a $U(1)$ bundle   satisfying:
\begin{enumerate}
\item There is a {\em bundle gerbe multiplication} which is a smooth isomorphism
$$
m \colon \pi_3^*(P) \otimes \pi_1^*(P) \to \pi_2^*(P)
$$
of $U(1)$ bundles over $Y^{[3]}$.    
  \item This multiplication is associative,   that is, if we let  $P_{(y_1, y_2)}$ denote the fibre of $P$ over $(y_1, y_2)$ then    the 
following diagram commutes for all $(y_1, y_2, y_3, y_4) \in Y^{[4]}$: 
\begin{equation*}
\begin{array}{ccc}
P_{(y_1, y_2) } \otimes P_{(y_2, y_3) } \otimes P_{(y_3, y_4) }  & \rightarrow  & P_{(y_1, y_3) } \otimes P_{(y_3, y_4) }  \\
\downarrow    &     & \downarrow    \\
P_{(y_1, y_2) } \otimes P_{(y_2, y_4) }  & \rightarrow &  P_{(y_1, y_4) }\\
\end{array}
\end{equation*}

\end{enumerate}
\end{definition}
  
We remark that for any $(y_1, y_2, y_3) \in Y^{[3]}$ the bundle gerbe multiplication defines an isomorphism:
$$
m \colon P_{(y_1, y_2) } \otimes P_{(y_2, y_3) }  \to P_{(y_1, y_3) } 
$$
of $U(1)$ spaces.
 
We can show using the bundle gerbe multiplication that there are natural
 isomorphisms $P_{(y_1, y_2)} \cong P_{(y_2, y_1)}^*$ and $P_{(y, y)} \simeq Y^{[2]} \times U(1)$.

We can rephrase the existence and associativity of the bundle gerbe multiplication to an
equivalent pair of conditions in the following way. The bundle gerbe multiplication 
gives rise to a section $s$ of $\delta(P) \to Y^{[3]}$.  Moreover $\delta(s)$ is a section of $\delta(\delta(P)) \to Y^{[4]}$.  But 
$\delta(\delta(P)) $ is canonically trivial so it makes sense to ask that $\delta(s) = 1$. This is the condition of associativity.  The family 
of spaces $\{ Y^{[p]} \mid p=1, 2, \dots \}$ is an example of a simplicial space \cite{Dup} and  by comparing to \cite{BryMac}
we see that  a bundle gerbe is the same thing as a {\em simplicial
line bundle} over this particular simplicial space.

\begin{example}
\label{ex:HC1}
If we replace $Y$ in the definition by $Y_\cU$ for some open cover $\cU$ of $M$ we obtain the definition of 
gerbe given by Hitchin 
\citeyear{MR1876068} and by his student Chatterjee  \citeyear{Cha}.  This consists of choosing an open 
cover $\cU$ of $M$ and a family of $U(1)$ bundles
$P \colon U_\a \cap U_\b$ such that over triple overlaps we have sections
$$
s_{\a\b\c} \in \Gamma(U_\a \cap U_\b \cap U_\c \mid P_{\b\c} \otimes P_{\a\c}^* \otimes P_{\a\b})
$$
and we require that $\delta(s) = 1$ in the appropriate way.  
\end{example}

\begin{example} 
\label{ex:3sphere}
The simplest example of a line bundle is given by the clutching construction on  the two sphere $S^2$. If $U_0$ and $U_1$ are the 
open neighbourhoods of the north and south hemispheres we take the transition function $g \colon U_0 \cap U_1  \to U(1)$ 
to have winding number one. As there are only two open sets there is no condition on triple overlaps and we obtain the 
$U(1)$ bundle over $S^2$ of chern class one.  In a similar fashion we can consider $U_0$ and $U_1$ to be 
open neighbourhoods of the north and south hemispheres of the three-sphere $S^3$. Their intersection is retractable 
to the two-sphere so we can choose over this the $U(1)$ bundle $P$ of chern class one. Again there are no additional conditions and
we obtain the gerbe of degree one over $S^3$. 
\end{example}

\begin{example}
\label{ex:HC2}
Hitchin and Chatterjee also consider a  gerbe as in Example \ref{ex:HC1} but with the
added requirement that  each $P_{\a\b} $
is trivial in the form $P_{\a\b} = U_\a \cap U_\b \times U(1)$. Writing elements of the disjoint union $Y^{[2]}_{\cU}$
as $(\alpha, \beta, x)$ where $x \in U_\a \cap U_\b$ we see that the bundle gerbe multiplication must take the form 
$$
((\alpha, \beta, x), z) \otimes ((\beta, \gamma, x), w) \mapsto ((\alpha, \gamma, x), zwg_{\a\b\c}(x) )
$$
for some $g_{\a\b\c} \colon U_\a \cap U_\b \cap U_\c \to U(1)$ and will be associative precisely when $g_{\a\b\c}$ is 
a co-cycle.  
\end{example}

We will refer to gerbes of the forms in Examples \ref{ex:HC1} or \ref{ex:HC2} as {\em Hitchin-Chatterjee gerbes}. The connection 
with bundle gerbes is simple.  For clarity we define:
\begin{definition} 
A bundle gerbe $(P, Y)$ over $M$ is called {\em local} if $Y = Y_\cU$ for some open cover $\cU$ 
of $M$.
\end{definition}
We then obviously have:
\begin{proposition}
A Hitchin-Chatterjee gerbe is the same thing as a local bundle gerbe.
\end{proposition}
 
If $(P, Y)$ is a bundle gerbe over $M$ then  associated to every point $m$ of $M$ we have a {\em groupoid} constructed as follows.
  The objects are the elements of the fibre $Y_m$ and the morphisms between $y_1$ and $y_2$ in $Y_m$ are  $P_{(y_1, y_2)}$.  
    Composition comes from the bundle gerbe multiplication.   If we call a groupoid a $U(1)$ groupoid if it is transitive and 
    the group of  morphisms of a point is 
isomorphic to $U(1)$,  then  the algebraic conditions on the bundle gerbe (that is the multiplication and its associativity)
 are captured precisely by saying that a bundle gerbe  is a bundle of $U(1)$ groupoids over $M$.

 We now consider the properties given in Section \ref{sec:intro} which we would like a $1$-gerbe to satisfy and show how 
 they are satisfied by bundle gerbes.  

  \subsection{Pullback}
If $f \colon N \to M$ then we can pullback $Y \to M$ to $f^*(Y) \to N$ with a map $\hat f \colon f^*(Y) \to Y$ covering $f$.   There
is an induced map ${\hat f}^{[2]} \colon  {f^*(Y)}^{[2]} \to Y^{[2]}$. Let 
$$
f^*(P, Y) =  ({{\hat f}^{[2]}}{}^*(P), f^*(Y)).
$$
To see this is a bundle gerbe notice that all this  is doing is pulling back the $U(1)$ groupoid at $f(n) \in M$ and
 placing it at $n \in N$ so we have a  bundle
  of $U(1)$ groupoids over $N$ and thus a bundle gerbe. 

\subsection{Dual and product}
If $(P, Y)$ is a bundle gerbe then $(P, Y)^* = (P^*, Y)$ is also
a bundle gerbe called the {\em dual} of $(P, Y)$.

  If $(P, Y)$ and $(Q, X)$ are bundle gerbes we can form the fibre product $Y \times_M X \to M$, a new surjective submersion 
    and then define a $U(1)$ bundle 
$$
P \otimes Q \to (Y \times_M X)^{[2]} = Y^{[2]} \times_M X^{[2]}
$$
  by 
$$
(P\otimes Q)_{((y_1, x_1), (y_2, x_2))} = P_{(y_1, y_2)} \otimes Q_{(x_1, x_2)}.
$$
  We define $(P, Y) \otimes (Q, X) = (P\otimes Q, Y\times_M X)$. 

\subsection{Characteristic class}
\label{subsec:DD}
The characteristic class of a bundle gerbe is called the {\em Dixmier-\-Douady}
class.      We construct it as follows.     Choose a good cover $\cU$ of $M$  \cite{BotTu} with 
sections $s_\a \colon U_\a \to Y$.      Then 
$$
(s_\a, s_\b) \colon U_\a \cap U_\b \to Y^{[2]}
$$
is a section.      Choose a section $\sigma_{\a\b} $ of $P_{\a\b} = (s_\a, s_\b)^*(P)$. That is some 
$$
\sigma_{\a\b} \colon  U_\a \cap U_\b \to P
$$
such that $\sigma_{\a\b}(x) \in P_{(s_\a(x), s_\b(x))}$.      Over triple overlaps we have 
$$
m(\sigma_{\a\b}(x), \sigma_{\b\c}(x)) = g_{\a\b\c}(x) \sigma_{\a\c}(x) \in P_{(s_\a(x), s_\c(x))}
$$
for  $g_{\a\b\c} \colon U_\a \cap U_\b \cap U_\c \to U(1)$.       This defines a 
co-cycle  which is the Dixmier-Douady class
$$
DD((P, Y)) = [g_{\a\b\c}] \in H^2(M, U(1)) = H^3(M, \ZZ).
$$

\begin{example}
If $\cU$ is an open cover and $g_{\a\b\c}$ a $U(1)$ co-cycle then we can build a Hitchin-Chatterjee gerbe 
or local bundle gerbe of the type considered in Example \ref{ex:HC2}.
It is easy to see that this has Dixmier-Douady class given by the \v Cech class $[g_{\a\b\c}]$.
\end{example}

Notice that this example shows that every class in $H^3(M, \ZZ)$ arises as the Dixmier-Douady class of some 
Hitchin-Chatterjee gerbe or of some (local) bundle gerbe. 

It is straightforward to check that if $f \colon N \to M$ and $(P, Y) $ is a bundle gerbe over $M$ then 
$f^*(DD(P, Y)) = DD(f^*(P, Y))$. Moreover we have

\begin{enumerate}
\item $DD((P, Y)^*)  = -DD((P, Y))$, and 
\item $DD((P, Y) \otimes (Q, X) ) = DD((P, Y)) + DD((Q, X))$.
\end{enumerate}

We will defer the question of triviality of a bundle gerbe until the next section and consider next the 
notion of a connective structure on a bundle gerbe. 

\subsection{Connective structure}
\label{sec:conn}
As  $P \to Y^{[2]}$ is a
$U(1)$ bundle  we can pick a connection $A$.     Call it a {\em bundle
gerbe connection } if it respects the bundle gerbe multiplication.   That is if the section $s$ of 
$\delta(P) \to Y^{[3]}$ satisfies $s^*(\delta(A)) = 0$, i.e is flat for $\delta(A)$. We would
like bundle gerbe connections to exist. This is a straightforward consequence of the  fact that the 
fundemental complex is exact.  Indeed if $A$ is any connection consider $s^*(\delta(A))$; we have
$\delta( s^*(\delta(A))) = \delta(s)^*(\delta\delta(A)) = 0$ because $\delta\delta(A)$ is the flat connection 
on the canonically trivial bundle $\delta\delta(P)$. Hence there is a one-form $a$ on $Y^{[2]}$ such that
$\delta(a) = s^*(\delta(A))$ and thus $A - a$ is a bundle gerbe connection. 

If $A$ is a bundle gerbe connection then the curvature $F_A \in \Omega^2(Y^{[2]})$ satisfies $\delta(F_A) = 0$.   
   From the exactness of the 
fundamental complex there must be an $f \in \Omega^2(Y)$ such
that $F_A = \delta(f)$.    As  $\delta$ commutes with $d$ we have  $\delta(df) = d\delta(f) = 
\allowbreak dF_A  \allowbreak = 0$.      Hence $df = \pi^*(\omega)$ for some $\omega \in \Omega^3(M)$.    So  
$\pi^*(d\omega) = \allowbreak d\pi^*(\omega) \allowbreak = d d f = 0$ and  $\omega $ is closed.      In fact it is a 
consequence of standard \v Cech de Rham theory that:
$$
\left[ \frac{1}{2\pi i} \omega\right] = r(DD((P, Y)))  \in H^3(M, \RR).
$$
    
We call $f$ a {\em curving} for $A$, the pair $(A, f)$ a 
{\em connective structure} for $(P, Y)$ and  $\omega$ is called the 
{\em three-curvature} of the connective structure $(A, f)$.  In string theory the 
two-form $f$ is called the $B$-field. 

We can give a local description of the connective structure as follows.  Assume we have an open 
cover $\cU$ of $M$ with local sections $s_\a \colon U_\a \to Y$ and sections over double overlaps
$\sigma_{\a\b}$ of $(s_\a, s_\b)^*(P) \to U_\a \cap U_\b$.  We define 
$$
A_{\a\b} = (s_\a, s_\b)^*(A) \in \Omega^1(U_\a \cap U_\b)
$$
and 
$$
f_\a = s_\a^*(f) \in \Omega^2(U_\a).
$$
These satisfy
\begin{align*}
A_{\b\c} - A_{\a\c} + A_{\a\b} &= g_{\a\b\c}^{-1} d g_{\a\b\c}\\
f_\b - f_\a &= df_{\a\b} \\
\end{align*}
and the three-curvature $\omega$ restricted to $U_\a$ is $df_\a$. 

\begin{example}
We can use this local description of the connective structure to calculate the Dixmier-Douady class of the Hitchin-Chatterjee  gerbe 
on the three sphere defined in Example \ref{ex:3sphere}. 
Stereographic projection from either pole identifies $S^3 - \{(1,0,0)\}$ and $S^3 - \{(-1,0,0)\}$ with $\RR^3$ and 
maps the equator to the unit sphere $S^2 \subset \RR^3$.  Let   $U_0$ and $U_1$ be  the pre-images of the interior
 of a ball of radius two in $\RR^3$ under both stereographic projections. 
We can identify  $U_0 \cap U_1$ with $S^2 \times (-1,1)$.  Pull back the 
line bundle of chern class $k$ on $S^2$, with connection $A$ and curvature $F$, to $U_0 \cap U_1$. Because there are  no 
triple overlaps this is a bundle gerbe connection. If we choose a partition of unity $\psi_0$ and $\psi_1$ for $U_0$ and 
$U_1$ then $f_0 = -\psi_1 F$  and $f_1 = \psi_0 F$ define two-forms on $U_0$ and $U_1$ respectively satisfying
$F = f_1 - f_0$ on $U_0 \cap U_1$. These two forms define a curving for the bundle gerbe connection.  The curvature is the globally defined 
three-form $\omega$ whose restriction to $U_0$ and $U_1$  is $-d\psi_1 \wedge F$ and $d\psi_0 \wedge F$ respectively.
The integral of $\omega$ over the three sphere reduces, by Stokes theorem, to the integral of $F$ over the two-sphere. 
Hence this bundle gerbe has Dixmier-Douady class $k \in H^3(M, \ZZ) = \ZZ$. 
\end{example}

Holonomy will need to wait until we have a considered the notion of triviality which we turn to now.
  
\section{Triviality} 

Recall that a $U(1)$ bundle  $P \to M$ is trivial if it is isomorphic to the bundle $M \times U(1)$ or, 
equivalently has a global section. This occurs if and only if $P \to M$ has zero Chern class.  If $s_a \colon U_\a \to P$ 
are local sections then $P$ is determined by a  transition function $g \colon U_\a \cap U_\b \to U(1)$ given by 
$s_\a = s_\b g_{\a\b}$ and $P \to M$ is trivial if and only if there exist $h_\a \colon U_\a \to U(1) $ such that 
$$
g_{\a\b} = h_{\b} h_{\a}^{-1}.
$$
In an analogous way  Hitchin and Chatterjee  \citeyear{Cha} define a gerbe $P_{\a\b} \to U_\a \cap U_\b$ to be   
 trivial if there are $U(1)$  bundles
$R_\a \to U_\a$ and isomorphisms $\phi_{\a\b} \colon R_\a \otimes R_\b^* \to P_{\a\b}$ on double-overlaps in such  a way
that the multiplication becomes the obvious contraction 
$$
R_\a \otimes R_\b^* \otimes R_\b \otimes R_\gamma^* \to R_\a \otimes R_\gamma^*.
$$

In the bundle gerbe formalism this idea takes the following form \cite{MurSte00Bundle-gerbes:-stable}. 
Let $R \to  Y$ be a $U(1)$ bundle and let $\delta(R) \to Y^{[2]}$ be defined as above.
Note that $\delta(R)$ has a natural associative bundle gerbe multiplication given by 
$$
\delta(R)_{(y_1, y_2)} \otimes \delta(R)_{(y_2, y_3)} = 
R_{y_1} \otimes R^*_{y_2} \otimes R_{y_2} \otimes R^*_{y_3} \simeq R_{y_1}  \otimes R^*_{y_3} = \delta(R)_{(y_1, y_3)}.
$$

\begin{definition} 
A bundle gerbe $(P, Y)$ over $M$ is called {\em trivial} if there is a $U(1)$  bundle $R \to Y$ such 
that $(P, Y)$ is isomorphic to $(\delta(R), Y)$.  In such a case we call a choice of $R$ and the
 isomorphism $\delta(R) \simeq P$ a {\em trivialisation} of $(P, Y)$. 
\end{definition}

\begin{example}
Let $(P, Y)$ be a bundle gerbe and assume that $Y \to M$ admits a global section $s \colon M \to Y$. 
Define $R \to Y $ by $R_y = P_{( s(\pi(y), y)}$. Then we have an isomorphism
\begin{align*}
\delta(R)_{(y_1, y_2)} &= P_{(s(\pi(y_2), y_2)} \otimes P^*_{(s(\pi(y_1), y_1)}\\
      &= P_{(s(\pi(y_2), y_2 )} \otimes P_{(y_1, s(\pi(y_1))}\\
      & \simeq P_{(y_1, y_2)}
\end{align*}
using the bundle gerbe multiplication and the fact that $s(\pi(y_1)) = s(\pi(y_2))$.  It is easy
to check that this isomorphism preserves the respective bundle gerbe multiplications and we have shown 
that if $Y$ admits a global section then any bundle gerbe $(P, Y)$ is trivial. Notice that the converse
 is not true. Just take an open cover with more than one element and $g_{\a\b\c} = 1$ to obtain a 
  Hitchin-Chatterjee gerbe which has zero Dixmier-Douady class but for which $Y_{\cU} \to M$ has no global section.
\end{example}

Consider the Dixmier-Douady class of $\delta(R)$. If $s_\a \colon U_\a \to Y$ are local sections for a good cover choose
local sections $\eta_\a $ of $s_\a^*(R)$. Then we can take as local sections of $(s_\a, s_\b)^*(\delta(R))$ the sections
$\sigma_{\a\b} = \eta_\a \otimes \eta_\b^* $ and it follows that the corresponding $g_{\a\b\c} = 1$ and $\delta(R)$ has
Dixmier-Douady class equal to zero.  The converse is also true. Consider a bundle gerbe with Dixmier-Douady class zero.
 So we have an open cover and $\sigma_{\a\b}$ such that 
$$
g_{\a\b\c} = h_{\b\c} h^{-1}_{\a\c} h_{\a\b}.
$$
By replacing $\sigma_{\a\b}$ by $\sigma_{\a\b}/h_{\a\b}$ we can assume that $g_{\a\b\c} = 1$.  Let $Y_\a = \pi^{-1}(U_\a)$ and define
$R_\a \to Y_\a$ by letting the fibre of $R_\a$ over $y \in Y_\a $ be $P_{(y,s_\alpha(\pi(y)))}$.  Construct an isomorphism
$\chi_{\a\b}(y)$ from the fibre of $R_\a$ over $y$ to the fibre of $R_\b$ over $y$ by noting that 
$$
P_{(y,s_\alpha(\pi(y)))} =   P_{(y,s_\beta(\pi(y)))}\otimes P_{(s_\alpha(\pi(y))s_\beta(\pi(y)))}
$$
and using $\sigma_{\a\b}(\pi(y)) \in P_{(s_\alpha(\pi(y))s_\beta(\pi(y)))}$.  Because $\sigma_{\b\c} \sigma_{\a\c} \sigma_{\a\b} = 1$
we can show that $\chi_{\a\b}(y) \circ \chi_{\b\c}(y) = \chi_{\a\c}(y)$ and  hence the $R_\a$ clutch together to form a global
$U(1)$  bundle $R \to Y$. It is straightforward to check that $\delta(R) = P$. 

Consider now a $U(1)$ bundle $R \to Y$ and assume that $\delta(R) \to Y^{[2]}$ has a section 
$s$ with $\delta(s) = 1$ with respect to the canonical trivialisation of $\delta(\delta(R)) \to Y^{[3]}$. 
The section $s$ is called {\em descent data} for $R$ and is equivalent to $R $ being the pull-back 
of a $U(1)$  bundle on $M$. Indeed $s$ constitutes a family of isomorphisms
$$
s(y_1, y_2) \colon R_{y_1} \to R_{y_2}
$$
and $\delta(s) = 1$ is equivalent to $s(y_2, y_3) \circ s(y_1, y_2) = s(y_1, y_3)$ from which it is easy 
to define a bundle on $M$ whose pull-back is $R$. 

Assume that a bundle gerbe $(P, Y)$ over $M$ is trivial and that $R_1$ and $R_2$ are two 
trivialisations. Then we have $\delta(R_1) \simeq P \simeq \delta(R_2)$ and hence a section 
of $\delta(R_1^* \otimes R_2)$ which is descent data for $R_1^* \otimes R_2$. Thus 
$R_1 = R_2 \otimes \pi^*(Q)$ for some $U(1)$ bundle $Q \to M$.  It is easy to show the converse that if $R$ is a trivialisation 
and $Q \to M$ a $U(1)$ bundle then $R \otimes \pi^*(Q)$ is another trivialisation. 
We have now proved

\begin{proposition}
Let $(P, Y)$ be a bundle gerbe over $M$. Then:
\begin{enumerate}
\item $(P, Y)$ is trivial  if and only if  $DD(P, Y) = 0$, and 
\item if $DD(P, Y) = 0$ then any two trivialisations of $(P, Y)$ differ by a $U(1)$ bundle 
on $M$.
\end{enumerate}
\end{proposition}

This should be compared to the 
case of $U(1)$ bundles ($0$-gerbes) where two trivialisations or sections of the bundle differ by a map into $U(1)$
which is a $-1$-gerbe. The general pattern is that we expect two trivialisations of a $p$-gerbe to differ by a $p-1$
gerbe.  Notice also that whereas any  two trivial $U(1)$ bundles are isomorphic there are many trivial bundle gerbes
which are not isomorphic. This leads us to the notion of stable isomorphism.

\begin{definition} 
If $(P, Y)$ and $(Q, X)$ are bundle gerbes over $M$ we say they are {\em stably 
isomorphic} \cite{MurSte00Bundle-gerbes:-stable}   if $(P, Y)^*\otimes(Q, X)$ is trivial. 
A choice of a trivialisation is called a {\em stable
 isomorphism} from $(P, Y)$ to $(Q, X)$. 
 \end{definition}

We have:

\begin{proposition} 
Bundle gerbes $(P, Y)$ and $(Q, X)$  over $M$ are stably isomorphic if and only if $DD(P, Y) = DD(Q, X)$.
The Dixmier-Douady class defines a bijection between stable isomorphism classes of bundle gerbes on $M$ and 
$H^3(M, \ZZ)$.
\end{proposition}
\begin{proof}
Bundle gerbes $(P, Y)$ and $(Q, X)$  over $M$ are stably isomorphic if and only 
if $(P, Y)^*\otimes(Q, X)$ is trivial which occurs if and only if $ - DD(P, Y) + D(Q, X) = 0$.  We have
 already seen that every 
three class arises as the Dixmier-Douady class of some bundle gerbe on $M$.
\end{proof}

It follows that the correct notion of equivalence for bundle gerbes is stable isomorphism. 
It is actually possible to compose stable isomorphisms and the details are given in work of
 Stevenson \citeyear{Ste} where the  structure of  the 
two category of all bundles gerbes on $M$ is  discussed.  See also  
\cite{Wal}. 

Note that we also have:
\begin{proposition}
Every bundle gerbe is stably isomorphic to a   Hitchin-Chatterjee gerbe.
\end{proposition}

\subsection{Holonomy}
\label{sec:holonomy}
Consider now a bundle gerbe $(P, Y)$ with connective structure $(A, f)$  over a surface $\Sigma$. 
Because $H^3(\Sigma, \ZZ) = 0$ we know that $(P, Y)$ is trivial. So there is a $U(1)$ bundle $R \to Y$ with 
$\delta(R) = P$. Choose a connection $a$ for $R$ and note that  $\delta(a)$ is a connection for $P$
 using the isomorphism $\delta(R) = P$. But $\delta(\delta(a)) $ is flat so $\delta(a)$ is a
  bundle gerbe connection. Hence $A = \delta(a) + \alpha$ for a one-form $\alpha$ on $Y^{[2]}$ with $\delta(\alpha) = 0$. Using the 
exactness of the fundamental complex we can solve $\alpha = \delta(\alpha')$ and hence show that 
$\delta(a + \alpha') = A$. So without loss of generality we can choose a connection $a$ on $R$ with 
$\delta(a) = A$.  Consider the two-form $ f - F_a$. This satisfies $\delta( f - F_a ) = F_A - F_{\delta(a)} = 0$
so $f - F_A = \pi^*(\mu_a)$ for some two form $\mu_a $ on $\Sigma$. Define the holonomy of $(A, f)$ over $\Sigma$ by 
$$
\hol((A, f), \Sigma) = \exp\left( \int_\Sigma \mu_a \right)
$$
and note that this is independent of the choice of trivialisation $R$ and connection $a$. Indeed any two trivialisations
with connection will differ by a $U(1)$ bundle on $M$ with connection and the corresponding $\mu_a$ will differ by the  curvature 
of the connection on that $U(1)$  bundle. But the integral of the curvature of a $U(1)$ bundle over a closed
surface is in $2\pi i \ZZ$ so the two definitions of holonomy agree.

If $(P, Y)$ is a bundle gerbe with connective structure $(A, f)$ on a general manifold $M$ and $\Sigma \subset M$ is 
a submanifold we can pull $(P, Y)$ and  $(A, f)$ back to $\Sigma$ and define $\hol((A, f), \Sigma)$ as above. In this 
more general setting if $X \subset M$ is a three dimensional submanifold with boundary $\partial X $ also a 
submanifold of $M$ we can trivialise $(P, Y)$ over all of $X$ and repeat the construction above. We then have
$d\mu_a = \omega$, the three-curvature of $(A, f)$, and thus
$$
\hol((A, f), \partial X) = \exp\left(\int_X \omega \right)
$$
the final property in Section \ref{sec:intro} which we wanted a $1$-gerbe to satisfy.

Assume that we have a local description for $(P, Y)$ and $(A, f)$ as in Section \ref{sec:conn} in terms 
of $g_{\a\b\c}$, $A_{\a\b}$ and $f_\a$. Then there is a remarkable formula \cite{MR1945806} for the holonomy,  first proposed in 
\cite{Alv} and also \cite{Gaw} and subsequently derived by a number of authors,  which can be described as follows.
  Choose a triangulation $\triangle$ of $\Sigma$ and a map 
$\chi \colon \triangle \to I$ such that for any simplex $\sigma \in \triangle$ we have $\sigma \subset U_{\chi(\sigma)}$.
Write $\sigma^2$, $\sigma^1$ and $\sigma^0$ for two, one and zero dimensional simplices, that is, faces, edges and vertices. Then

\begin{multline*}
\hol((A, f), \Sigma) = \\
\exp\left( \sum_{\sigma^2 } \int_{\sigma^2} f_{\chi(\sigma^2)}\right) \exp\left(
\sum_{\sigma^1 \subset \sigma^2} \int_{\sigma^1} A_{\chi(\sigma^2)\chi(\sigma^1)} \right) \prod_{\sigma^0 \subset \sigma^1 \subset \sigma^2} g_{\chi(\sigma^2)\chi(\sigma^1)\chi(\sigma^0)}(\sigma^0).
\end{multline*}

\subsection{Obstructions to certain kinds of $Y \to M$.}
When we consider the examples in the next section it will become apparent that they tend to cluster into 
two kinds:  either $Y \to M$ has infinite dimensional fibres or $Y \to M$ has discrete fibres as in the
 case of Hitchin-Chatterjee gerbes. There
is a reason for this which is the following result.

\begin{proposition} 
\label{prop:nogo}
Let $(P, Y)$ be a bundle gerbe over $M$ with $Y \to M$ a finite
dimensional fibration with both $M$ and the fibres of $Y \to M$ one connected. Then the Dixmier-Douady class of $(P, Y)$ is torsion.
\end{proposition}
\begin{proof} 
The proof uses the result from \cite{Got} which shows that if $Y \to M$ satisfies the
 hypothesis we have stated in the Proposition and moreover there is a smooth two-form $\mu$  defined
  on the vertical tangent bundle of $Y \to M$ which is closed on each fibre then we can extend $\mu$  to a 
  global, closed two form on $Y$.  To apply this choose a 
connection and curving $f$ for $(P, Y)$. If we restrict $f$ to any fibre
of $Y \to M$ it agrees with $F$ and  hence is closed in the fibre directions. So there
exists a global two form $\mu$ on $Y \to M$ which agrees with $f$ in the vertical directions. 
Consider $\rho = f - \mu$. Both $\rho$ and $d\rho$ are zero in the vertical 
directions so that $\rho = \pi^*(\chi)$ for some two-form $\chi$ on $M$. But then 
$\pi^*(d\rho) = d(f - \mu) = \pi^*(\omega)$ where $\omega $ is the three curvature. Hence the 
image of the Dixmier-Douady class in real cohomology is zero so the Dixmier-Douady class is torsion. 
\end{proof}

\section{Examples of bundle gerbes}
We expect to find bundle gerbes on manifolds $M$ which have three dimensional 
cohomology. There are two, related, cases we will discuss here. The first comes from  
lifting problems for  principal bundles, the so-called {\em lifting bundle gerbe} and the second is to take $M$ a simple,  compact Lie group. 

\subsection{Lifting bundle gerbe}
Consider a central extension of Lie groups
$$
0 \to U(1) \to \hat\cG \stackrel{\rho}{\to} \cG \to 0.
$$
so that $U(1)$ is the kernel of $p$ and is in the center of $\cG$. 
If $Y \to M$ is a $\cG$ principal bundle then we can ask if it lifts to a $\hat \cG$
bundle $\hat Y \to M$.  That is can we find a $\hat\cG$ bundle $\hat Y \to M$  and a bundle
 morphism $f \colon \hat Y \to Y$ such that 
$f(pg) =  f(p) \rho(g)$ for all $p \in \hat Y$ and $g \in \hat G$.  There is a  well known topological
obstruction to the existence of $\hat Y$ which we can calculate as follows. Choose a good open cover $\cU$ of 
$M$ and local sections of $Y$ that give rise to a transition function 
$$
g_{\a\b}  \colon U_\a \cap U_\b \to G
$$
in the usual way. Because these double overlaps are all contractible we can choose 
lifts of each $g_{\a\b}$ to $\hat g_{\a\b} \colon U_\a \cap U_\b \to \hat G$.  Note
that 
$$
\epsilon_{\a\b\c} = g_{\b\c} g^{-1}_{\a\c} g_{\a\b} 
$$
takes values in $U(1)$. It is not difficult to show that $\epsilon_{\a\b\c}$ is a 
$U(1)$ valued co-cycle and that the class
$$
[\epsilon_{\a\b\c}] \in H^2(M, U(1) \simeq H^3(M, \ZZ) 
$$
vanishes if and only if $P$ lifts to a $\hat \cG $ bundle.

Given a $\cG$ bundle $Y \to M$ there is a map $g \colon Y^{[2]} \to \cG$ defined by $y_1 g(y_1, y_2) = y_2$.  
Notice that $\hat\cG \to \cG$ is a $U(1)$ bundle so we can pull it back 
to define a $U(1)$ bundle $Q \to Y^{[2]} $ whose fibres are cosets of $U(1)$ in $\hat G$ defined by 
$$
Q_{(y_1, y_2)} = U(1)g(y_1, y_2)
$$
Because  $g(y_1, y_2) g(y_2, y_3) =  g(y_1, y_3)$ the product of an element in the coset containing $g(y_1, y_2)$ 
with an element in the coset containing $g(y_2, y_3)$ will be an element in the coset containing $g(y_1, y_3)$ which
defines a bundle gerbe multiplication 
$$
Q_{(y_1, y_2)} \otimes Q_{(y_2, y_3)} \to Q_{(y_1, y_3)}.
$$
The bundle gerbe $(Q, Y)$ is called the {\em lifting bundle gerbe} of $Y$ \cite{Mur96Bundle-gerbes}.
 It is easy to check that the lifting  bundle gerbe has Dixmier-Douady class precisely the 
obstruction to lifting the bundle $Y \to M$. Indeed if we follow through the construction in Section \ref{subsec:DD} we
find that $\hat \sigma_{\a\b} = \hat g_{\a\b}$.  It follows from the discussion above that the lifting bundle
gerbe is trivial if and only if the bundle $Y \to M$ lifts to $\hat \cG$. In fact this follows directly because a lift 
$\hat Y \to Y$ will be a $U(1)$ bundle over $Y$ and actually a  trivialization of the lifting bundle gerbe defined above.

\subsection{Projective bundles}
\label{sec:proj}
Let $H$ be a Hilbert space, possibly finite dimensional. A Hilbert  bundle with fibre $H$ can be
 regarded locally as a collection of transition functions 
$$
g_{\a\b} \colon U_\a \cap U_\b \to U(H)
$$
satisfying the co-cycle condition
$$
g_{\b\c} g_{\a\c}^{-1} g_{\a\b} = 1.
$$
In some situations we have slightly less than this, namely a collection of transition functions 
$$
g_{\a\b} \colon U_\a \cap U_\b \to U(H) 
$$
satisfying 
$$
g_{\b\c} g_{\a\c}^{-1} g_{\a\b} =  \epsilon_{\a\b\c} 1_{U(H)}.
$$
If we denote by $\rho \colon U(H) \to PU(H)$ the projection onto the projective unitary group   we see that
 have
$$
\rho(g_{\b\c}) \rho(g_{\a\c}^{-1}) \rho(g_{\a\b}) =  \rho(\epsilon_{\a\b\c} 1) = 1_{PU(H)}
$$
so there is a  well-defined bundle of projective spaces or a principal $PU(H)$ bundle. Given a
 projective bundle a natural question is to ask when it is the projectivisation of a global Hilbert bundle. This is equivalent to lifting
the $PU(H)$ bundle to $U(H)$. The obstruction to this lifting is the class is  
$$
[\epsilon_{\a\b\c}] \in H^2(M, U(1) )
$$
arising from the sequence
$$
0 \to U(1) \to U(H) \to PU(H) \to 0.
$$

In the finite dimensional class we can also consider the obstruction to lifting from 
$PU(n)$ to $SU(n)$ via the exact sequence
$$
0 \to \ZZ_n \to SU(n) \to PU(n) \to 0
$$
In this case  the lifting bundle gerbe is really a $\ZZ_n \subset U(1)$ bundle gerbe 
and the Dixmier-Douady class is a torsion class in the image of the Bockstein map
$$
H^2(M, \ZZ_n ) \to H^3(M, \ZZ).
$$
For further details on $\ZZ_n$ bundle gerbes see \cite{CarMicMur00Bundle-gerbes-applied}. Note finally that 
if $Y \to M$ is a $PU(n)$ bundle then the lifting bundle gerbe has finite dimensional fibres so the fact that
its Dixmier-Douady class is torsion is implied by  Proposition \ref{prop:nogo}.

\subsection{Bundle gerbes on Lie groups}
If $G$ is a compact, simple,  Lie group then $H^3(G, \ZZ) \simeq \ZZ$ so we expect to find bundle gerbes 
on $G$. There are, in fact,  a number of constructions of bundle gerbes on $G$ and we have seen one already
 for $SU(2) \simeq S^3$ in Section \ref{ex:3sphere}.  For convenience let us fix an orientation of $G$ and call the bundle gerbe over $G$ 
of Dixmier-Douady class one the {\em basic bundle gerbe} on $G$.

\begin{example}
Let $PG$ be the based path space of all smooth maps $g$ of the interval 
$[0, 1]$ into $G$ with $g(0) = 1$. Let $\ev \colon PG \to G$ be the 
evaluation map $\ev(g) = g(1)$. The kernel of $\ev$ is the space of based
loops $\Omega(G)$, that is all smooth maps $g \colon [0, 1] \to G$ 
for which $g(0) = 1 = g(1)$.  In the smooth Fr\'echet topology 
$PG \to G$ is a principal $\Omega(G)$ bundle. Moreover there is a well-known 
central extension \cite{PreSeg}
$$
0 \to U(1) \to \widehat \Omega(G) \to \Omega(G) \to 0
$$
where $\Omega(G) $ is the Kac-Moody group. Hence there is a corresponding 
lifting bundle gerbe over $G$ which is the basic bundle gerbe. It is possible to give an explicit construction 
of this gerbe over $G$, see \cite{Mur96Bundle-gerbes} and also \cite{Ste}.
\end{example}

If we wish to avoid the infinite dimensional spaces there is a  construction 
of the basic bundle gerbe over $G$ in the simply connected case due to 
Meinrenken \citeyear{Men} and see also \cite{MR2078218} for the non-simply connected case.  
This construction has disconnected fibres for $Y \to G$ and uses the standard structure theory of compact, simple, simply-connected Lie groups.

\begin{example}
\label{ex:unitary}
For $G = SU(n)$ there is a simple construction of the basic bundle gerbe due to Meinrenken \citeyear{Men} (see also \cite{Mic}) using an
 open cover of $G$ which can be presented without the 
cover as follows. Let
$$
Y = \{ (g, z) \mid \det(g - z1) \neq 0 \} \subset G \times U(1).
$$
For convenience let us write an element $((g, z), (g, w)) \in Y^{[2]}$ 
as $(g, w, z)$. If $u \in U(1)$ and $u \neq w$ and $u \neq z$ let us 
say that $u$ is between $w$ and $z$ if an anti-clockwise rotation of $z$ into 
$w$ passes through $u$.  Then let $W_{(g, w, z)}$ be the sum of all the
eigenspaces of $g$ for eigenvalues between $z$ and $w$ and define $P_{(g, w, z)}$ to be the $U(1)$ frame 
bundle of $\det W_{(g, w, z)}$. To define the bundle gerbe product notice that if $u$ is between $w$ and $z$ then 
$$
W_{(g, w, u)} \oplus W_{(g, u, z)}  \oplus W_{(g, z, w)} = \CC^n 
$$
so that 
$$
\det W_{(g, w, u)}  \otimes \det W_{(g, u, z)} \otimes \det W_{(g, z, w)} = \CC
$$
and thus
$$
\det W_{(g, w, u)}  \otimes \det W_{(g, u, z)}  =  \det W_{(g, z, w)}^* 
$$
Similarly 
$$
W_{(g, w, z)}  \oplus W_{(g, z, w)} = \CC^n 
$$
so that 
$$
\det W_{(g, w, z)}  \otimes \det W_{(g, z, w)} = \CC
$$
and 
$$
\det W_{(g, w, z)}  = \det W_{(g, z, w)}^*.
$$
There are a number of other cases that can be dealt with in a similar fashion and putting all these facts together 
gives a bundle gerbe multiplication on $P \to Y^{[2]}$.  A construction of the curving on $(P, Y)$ will appear in \cite{MurSte}.
\end{example}

\section{Applications of bundle gerbes}

\subsection{The Wess-Zumino-Witten term} 
The Wess-Zumino-Witten term associates to a smooth map $g$ of a surface $\Sigma$ into a compact, simple,  Lie 
group $G$ an invariant $\Gamma(g)
\in U(1)$. As noted by a number of authors \cite{CarMicMur00Bundle-gerbes-applied,MR1887879,MR1945806} this can be
 understood as the holonomy of a connection and curving on the  basic gerbe on the group.   The original definition
 of Witten \citeyear{Wit} is that we choose 
a three-manifold
$X$ with $\partial X = \Sigma$ and extend $g$ to $\hat g \colon X \to G$. We then consider 
$$
\int_X \hat g^*(\omega)
$$ 
where $\omega$ is a three form on $G$ representing a generator of $H^3(G, 2\pi i\ZZ)$.
If we choose a different extension $\tilde g \colon X \to G$ then the 
pair can be combined to define a map from the manifold $X \cup_\Sigma X$, formed
by joining two copies of $X$ (with opposite orientations) along $\Sigma$, into $G$. 
Call this map $\hat g \cup \tilde g$. Then
$$
\int_X \hat g^*(\omega) - \int_X \tilde g^*(\omega) = 
\int_{X \cup_\Sigma X} (\hat g \cup \tilde g)^*(\omega)  \in 2\pi i\ZZ.
$$
It follows that 
$$
\Gamma(g) = \exp\left( \int_X \hat g^*(\omega)\right) \in U(1)
$$
depends only on $g$. 

If we choose a suitable connection and curving $(A, f)$ on the basic gerbe on $G$, so that 
it has curvature $\omega$ then we see that 
$$
\Gamma(g)  = \hol(\Sigma, g^*(A, f)).
$$
The gerbe  approach to the Wess-Zumino-Witten term has two advantages. Firstly it removes the topological restriction 
on $M$ of $2$-connectedness necessary in Witten's definition so that the map $g$ can be extended to the three-manifold $X$. 
Secondly  we can use the local formula for the holonomy given in Section \ref{sec:holonomy}. For details see \cite{MR1945806}.

\subsection{The Faddeev-Mickelsson anomaly}
We follow  \cite{CarMur96Faddeevs-anomaly-and-bundle} and \cite{Seg}. Let $X$ be a compact, Riemannian,
spin, three-manifold and denote by $\cA$ the space of connections on a complex vector bundle over $M$ and
 by $\cG$ the space of gauge transformations.  For any $A \in \cA$ the chiral Dirac operator $D_A$ coupled to
  $A$ has discrete spectrum. Let 
$$
Y = \{ (A, t) \mid t \notin \spec(D_A) \}
$$
considered as a submersion over $\cA$. Note that $\cG$ acts on $Y$ and $Y/\cG \to \cA/\cG$ is another submersion. Following
Segal \citeyear{Seg} for $(A, s) \in Y$ we can decompose the Hilbert space $H$ of coupled spinors into a direct sum of 
eigenspaces of $D_A$ for eigenvalues greater than $s$ and a direct sum of eigenspaces of $D_A$ for eigenvalues less than $s$.
Denote these by $H_{(A, s)}^-$ and $H_{(A,s)}^+$ respectively.   We can then form the Fock space
$$
F_{(A, s)} = \bigwedge H_{(A,s)}^+ \otimes \bigwedge (H_{(A,s)}^-)^*
$$
which is a bundle $F \to Y$. Notice that $\cG$ acts on $F$ and gives rise to a bundle $F/\cG \to Y/\cG$. 

If we choose another $t < s$ we have
$$
H = H_{(A,t)}^- \oplus V_{(A, t, s)} \oplus H_{(A,s)}^+
$$
where $V_{(A, t, s)}$ is the sum of all the eigenspaces of $D_A$ for eigenvalues between $t$ and $s$. If we use the canonical
isomorphism 
$$
\bigwedge V^*_{(A, t, s)} \otimes \det V_{(A, t, s)} = \bigwedge V_{(A, t, s)}
$$
then we can show that 
$$
F_{(A, s)} = F_{(A, t)} \otimes \det V_{(A, t, s)}.
$$
It follows that the projective spaces of $F_{(A, s)}$ and $F_{(A, t)}$ are canonically isomorphic and descend to
 a projective bundle $\cP \to \cA$ and as $\cG $ acts there is also a projective bundle $\cP/ \cG \to \cA/\cG$. 
  The question of interest \cite{Seg} is whether there is a Hilbert bundle $\cH \to \cA/\cG$ whose projectivisation 
  is $\cP / \cG$. Notice that as $\cA$ is 
contractible the answer to the equivalent question on $\cA$ is clearly positive. 

As noted in Section \ref{sec:proj} we  could give a bundle gerbe interpretation of this question via the lifting
 bundle gerbe for the central extension
$$
0 \to U(1) \to  PU(\tilde \cH) \to U(\tilde \cH) \to 0
$$
for a suitable Hilbert space $\tilde \cH$. However there is a more direct approach as follows:
Let $P_{(A, s, t)}$ be the unitary frame bundle of $\det V_{(A, t, s)}$. Notice that if $r < s < t$ then 
$$
 V_{(A, r, t)}\oplus  V_{(A, t, s)} =  V_{(A, r, s)}
$$
so that 
$$
\det V_{(A, r, t)} \otimes  \det V_{(A, t, s)} =  \det V_{(A, r, s)}
$$
gives rise to a bundle gerbe multiplication similar to the $SU(n)$ case Example \ref{ex:unitary}.  Again $\cG$ 
acts so this descends to a bundle gerbe on $\cA/\cG$. If this bundle gerbe is trivial then there is a line
bundle $L \to Y$ such that 
$$
\det V_{(A, t, s)} = L_{(A, s)} \otimes L_{(A, t)}^*
$$
and moreover as it is the bundle gerbe on $\cA/\cG $ which is trivial these are $\cG$ equivariant isomorphisms. 
It follows that 
$$
F_{(A, s)} = F_{(A, t)} \otimes L_{(A, s)} \otimes L_{(A, t)}^*
$$
or 
$$
F_{(A, s)}\otimes L_{(A, s)}^* = F_{(A, t)}  \otimes L_{(A, t)}^*.
$$
and these are $\cG$ equivariant isomorphisms. Hence $F_{(A, s)}\otimes L_{(A, s)}^*$ descends to a 
$\cG$ equivariant Hilbert bundle on $\cA$ whose projectivization is $\cP$. It follows that there is a 
Hilbert bundle on $\cA/\cG$ whose projectivization is $\cP/\cG$. 

\subsection{String structures}
In \citeyear{Kil} Killingback introduced the notion of a {\em string structure}. Given a  $G$ bundle $Q \to X$ 
he considers the corresponding loop bundle $LQ \to LX$ which is an $LG$ bundle. As noted above
we have the Kac-Moody central extension 
$$ 
0 \to U(1) \to \widehat {LG} \to LG \to 0.
$$
Killingback says that $Q \to X$ is {\em string} if $LQ \to LX$ lifts to $\widehat {LG}$,  calls a choice of such a lift 
a {\em string structure} and defines the obstruction class to the lift to be the {\em string class}. 
 In \cite{MurSte03Higgs-fields-bundle} we consider the more general situation 
of an $LG$ bundle on a manifold $M$ and the question of whether it lifts to $\widehat{LG}$.  There is an equivalence
between $LG$ bundles on $M$ and $G$ bundles on $S^1 \times M$ which is exploited to define a connection and 
curving on the lifting bundle gerbe associated to the $LG$ bundle. This enables the derivation of a de Rham 
representative for the image of the string class in real cohomology.  The same approach could 
be applied to the other level Kac-Moody central extensions.

\section{Other matters}
In the interests of brevity nothing has been said about the original approach to gerbes as sheaves of groupoids. 
Details are given in \cite{Bry} and the relationship with bundle gerbes is discussed in \cite{Mur96Bundle-gerbes} and in more detail 
in \cite{Ste}.  In this same context I  would like to thank Larry Breen for pointing out that in the work of Ulbrich \citeyear{MR1049754} and \citeyear{MR1104600} the   notion of cocycle bitorsors can be interpreted as a form of a bundle gerbe.

In the definitions and theory above we could replace $U(1)$ by any abelian group $H$ and
 there would only be  obvious minor modifications such as the Dixmier-Douady class being 
 in $H^2(M, H)$. If we want to replace $H$ by a non-abelian group things become more difficult 
 as we mentioned in the introduction because we cannot form the product of two $H$ bundles if $H$ 
 is non-abelian. To get around this difficulty we need to replace principal bundles by principal bibundles 
 which have a left and right group action. The resulting theory becomes more complicated although closer to 
 Giraud's original aim of understanding non-abelian cohomology.
 For details see  \cite{MR2117631} and \cite{MR2183393}.

We have motivated gerbes by the idea of  replacing the  transition function $g_{\a\b} $
by a $U(1)$ bundle $P_{\a\b}$ on double overlaps. It is natural to consider what happens if we  take the next
 step and replace $P_{\a\b}$  by a bundle 
gerbe  on each double overlap. This gives rise to the notion of bundle two gerbes whose characteristic class 
is a four class. In particular there is associated to any principal $G$ bundle $P \to M$ a bundle two gerbe whose
 characteristic class is the 
pontrjagin class of $P$.  For details of the theory see \cite{Ste,MR2032513} and for an application
to Chern-Simons theory see \cite{CarJohMurWan05Bundle-gerbes-for-Chern-Simons}.

It is clear that we can continue on in this fashion and consider bundle two gerbes on double overlaps and 
more generally inductively use $p$ gerbes on double overlaps to define $p+1$ gerbes. However the theory becomes 
increasingly complex for a reason that we have not paid much attention to in the discussion above. In the case of $U(1)$ bundles
the transition function has to satisfy one condition, the co-cycle identity. In the case of gerbes the $U(1)$ bundles
over double overlaps satisfy two conditions, the existence of a bundle gerbe multiplication and its associativity. 
In the case of bundle two gerbes there are three conditions and the complexity continues to grow in this fashion.


\end{document}